\documentclass[12pt]{article}
\usepackage{amsmath,amssymb}
\usepackage{rotating}
\usepackage{xypic}
\usepackage{amssymb}
\usepackage{amsmath,amscd}
\usepackage{pst-node,pstricks,multido,pst-plot,pst-text,pst-3d}%
\usepackage{graphicx}
\usepackage{amsfonts}
\usepackage{authblk}
\usepackage{float}

\usepackage{appendix}
\usepackage{titlesec}
\titlespacing*{\section}{0pt}{1.1\baselineskip}{\baselineskip}
\usepackage{footmisc}
\newtheorem{example}{Example}[section]

\newtheorem{remark}[example]{Remark}
\newtheorem{theorem}[example]{Theorem}

\newtheorem{corollary}[example]{Corollary}

\newtheorem{proposition}[example]{Proposition}

\newtheorem{lemma}[example]{Lemma}
\usepackage[colorlinks=true,linktoc=all,linkcolor=blue]{hyperref}

\def\<{\langle}
\def\>{\rangle}

\def\N{{\mathbb N}}

\def\ashuff#1#2#3{
\kern 1pt \vrule height#1 \overline{\vrule height#3 width 0pt
\hskip#2} \rule{.3pt}{#1}\overline{\vrule height#3 width 0pt
\hskip#2} \rule{.3pt}{#1} \kern 1pt }

\makeatletter
\let\@fnsymbol\@arabic

\def\Carre3#1{\left[\begin{array}{ccc}#1\end{array}\right]}

\def\pdiff#1#2{\left(\partial #1\over\partial #2\right)}
\def\diff#1#2{{\partial #1\over\partial #2}}
\def\idiff#1#2#3{{\partial^{#3} #1\over\partial #2^{#3}}}
\begin{document}

\title{On the self-convolution of generalized Fibonacci numbers}
\author{Hac\`ene Belbachir\footnote{{\tt hacenebelbachir@gmail.com},\\ USTHB, Faculty of Mathematics, RECITS Laboratory  BP 32, El Alia 16111 Bab Ezzouar, ALGIERS, ALGERIA.}, Toufik Djellal\footnote{{\tt \text{t\_djellal@esi.dz}}, \\ USTHB, Faculty of Mathematics, RECITS Laboratory  BP 32, El Alia 16111 Bab Ezzouar,  ALGIERS, ALGERIA. \\ ESI, BP 108M Oued Smar, 16309, El Harrach, ALGIERS, ALGERIA.}, Jean-Gabriel Luque\footnote{{\tt jean-gabriel.luque@univ-rouen.fr}, \\ Normandie Universit\'e. Universit\'e de Rouen. Laboratoire LITIS  EA 4108. 
76800 SAINT-\'ETIENNE DU ROUVRAY FRANCE.}}

\maketitle 
\begin{abstract}
We focus on a family of equalities pioneered by Zhang  and generalized by Zao and Wang 
 and hence by Mansour  which involves self convolution of generalized Fibonacci numbers.
 We show that all these formulas are nicely stated in only one equation involving a bivariate ordinary generating function and we give also a formula
 for the coefficients appearing in that context. As a consequence, we give the general forms for the equalities of Zhang, Zao-Wang and Mansour.
 \end{abstract}
\section{Introduction}
Generalized Fibonacci numbers and generalized Lucas numbers commonly refer to two sequences (resp. $U_{k}$ and $V_{k}$) with the same 
characteristic polynomial $Q=x^{2}-px+q$ which has two distinct roots  such that $U_{0}=0$, $U_{1}=1$, $V_{0}=2$, and $V_{1}=p$,
 \cite{Horadam,belb4}. Under this assumption the sequences have the following Binet form
 \begin{equation}
	 U_{k}={\alpha^{k}-\beta^{k}\over \alpha-\beta},\,V_{k}=\alpha^{k}+\beta^{k}.
 \end{equation}
 Generalized Fibonacci numbers have several applications in numerical analysis for instance for solving nonlinear equations, 
 (see e.g. \cite{Jamieson,MP}).  One of the authors also investigated closed topics such as  sums of product of Fibonacci polynomials \cite{belb1}, 
q-analogs \cite{belb,belb2,belb3}.
 
  We focus on a family of equalities pioneered by Zhang \cite{Zhang} and generalized by Zao and Wang \cite{ZW}
 and hence by Mansour \cite{Mansour} which involves self convolution of Fibonacci numbers. The starting point of these equalities is that
 the discriminant of a quadratic form can be written as a differential operator.
 In particular, Mansour \cite{Mansour}
 gave  a method to obtained all the equations but did not provide a closed form. 
 We show that all these formulas are nicely stated in only one equation involving a bivariate ordinary generating function and we give also a formula
 for the coefficients appearing in that context. As a consequence, we give the general forms for the equalities of Zhang, Zao-Wang and Mansour.\\

In the last section, we explain how to generalize these results for higher level recurrence by the use of classical invariant theory and  we detail the example of Tribonacci numbers.
 Throughout the paper we use extensively materials issued from the theory of generating functions; readers
  should refers to \cite{Flajolet,Wilf} for a survey of the topic. 
\section{A polarization formula}
Consider a quadratic polynomial $Q$ in $x$ with two distinct roots. 
Without lost of generality we assume $Q$ unitary (i.e. ${\partial^{2}\over \partial x^{2}}Q=2$). If $\Delta$ denotes the discriminant of
$Q$, one has
\begin{equation}\label{Q2Delta}
	\Delta=\left({\partial\over\partial x }Q\right)^{2}-4Q.
\end{equation} As a consequence, setting $F=-\frac1Q$, one has
\begin{equation}\label{F^2}
	F^{2}=\frac1\Delta\left(\left({\partial\over\partial x }Q\right){\partial\over\partial x }F+4F\right).
\end{equation}
Our purpose is to describe the polynomials  $\alpha_{i,n}(x)$ such that
\begin{equation}\label{defalpha}
	F^{n}=\sum_{i=0}^{n-1}\alpha_{i,n}(x){\partial^{i}\over\partial x^{i}}F,
\end{equation}
which are obtained  applying successively many times (\ref{Q2Delta}).
We consider the generating function
\begin{equation}
	\mathcal F(y)=\sum_{n\geq 1}F^{n}y^{n}={yF\over 1-yF}.
\end{equation}
We want to find a differential operator $\mathbb D_{x}^{y}$  satisfying
\begin{equation}\label{defDxy}
	\mathcal F(y)=\mathbb D_{x}^{y}F.
\end{equation}
Remark first, that if we replace $x$ by $x=x+x'$ with $x'=\frac12\left(\frac1{\sqrt{1-{4y\over\Delta}}}-1\right){\partial Q\over\partial x}$ in $Q$,
 we obtain
\begin{equation}
	Q\left(x'\right)={Q+y\over 1-4{y\over\Delta}}.
\end{equation}
We translate this equality in terms of differential operators considering the polarization operator $\mathbb P_{X}$ sending each $f(x)$ to 
$f(x+X)$. This operator is
\begin{equation}
\mathbb P_{X}=\exp\left\{X{\partial\over\partial x}\right\}.
\end{equation}
We set 
\begin{equation}\label{Dxy}
\mathbb D_{x}^{y}=\frac y{1-4{y\over\Delta}}\mathbb P_{X}\left|_{X\rightarrow x'}\right.,
\end{equation}
where $\mathbb P_{X}\left|_{X\rightarrow x'}\right.$ means that we first act by the operator 
 $\mathbb P_{X}$ and hence we apply the substitution $X\rightarrow x'$.
This operation is necessary because $\partial Q\over\partial x$ is not a constant. So we obtain
\begin{equation}\label{D2F}
\mathbb D_{x}^{y}F
= -\frac y{Q+y}=
-\frac {\frac yQ}{1+\frac yQ}={yF\over 1-yF}=\mathcal F(y).
\end{equation}
Expanding $\mathbb D_{x}^{y}$ and comparing to equality (\ref{defalpha}), we deduce our main result:
\begin{theorem}\label{theoalpha}
	We have
	\begin{equation}	\alpha_{i,n}={\left(\partial Q\over\partial x\right)^{i}\over \Delta^{n-1}}\tilde\alpha_{i,n},
	\end{equation}
	where the ordinary generating function of the $\tilde\alpha_{i,n}$'s is
	\begin{equation}\label{A(y,z)}
		A(y,z)=\sum_{0\leq i<n}\tilde\alpha_{i,n}z^{i}y^{n}={y\over1-4y}\exp\left\{\frac12z\left(\frac1{\sqrt{1-4y}}-1\right)\right\}.
	\end{equation}
\end{theorem}

\begin{example}\rm{}
The first terms of the expansions of $A(y,z)$ give
\[{}\begin{array}{rcl}
A(y,z)&=&y+ \left( z+4 \right) {y}^{2}+ \left( 7\,z+\frac12\,{z}^{2}+16 \right) {y
}^{3}\\&&+ \left( 38\,z+5\,{z}^{2}+\frac16\,{z}^{3}+64 \right) {y}^{4}\\&&+
 \left( 187\,z+{\frac {69}{2}}\,{z}^{2}+{\frac {13}{6}}\,{z}^{3}+\frac1{24}
\,{z}^{4}+256 \right) {y}^{5}+\cdots .
\end{array}
\]
From Theorem \ref{theoalpha} we deduce
\begin{itemize}
	\item $\Delta F^{2}=\diff Qx\diff Fx+4F$ 
	\item $\Delta^{2}F^{3}=\frac12{\pdiff Qx}^{2}\idiff Fx2+7\diff Qx\diff Fx+16F$
	\item $\Delta^{3}F^{4}=\frac16{\pdiff Qx}^{3}\idiff Fx3+5{\pdiff Qx}^{2}\idiff Fx2
+38\diff Qx\diff Fx+64F$
\item $\Delta^{4}F^{5}=\frac1{24}{\pdiff Qx}^{4}\idiff Fx4+\frac{13}6{\pdiff Qx}^{3}\idiff Fx3
+\frac{69}2{\pdiff Qx}^{2}\idiff Fx2+187\diff Qx\diff Fx+256F$
\item $\dots$
\end{itemize}
\end{example}
\begin{remark}\rm
	In our proof of Theorem \ref{theoalpha}, the fact that Equality (\ref{A(y,z)}) is obtained by applying
	successively many times formula \ref{Q2Delta} is somewhat hidden.  In appendix \ref{AppProof} we give an alternative proof which
	makes this property more visible. This proof is based on a solving of a bivariate differential equation.
\end{remark}

\section{Itering self-convolution}
Consider two sequences of complex numbers $C=(c_{n})_{n\geq 0}$ and $D=(d_{n})_{n\geq 0}$. The convolution of $C$ and $D$
is the sequence defined by 
\begin{equation}C\star D=\left(\sum_{i+j=n}c_{i}d_{j}\right)_{n\in\N}.\end{equation} In other words, the ordinary generating 
function of the convolution of two sequences is the product of the ordinary generating function of the sequences,
\begin{equation}
	S_{C\star D}(z)=S_{C}(z)S_{D}(z).
\end{equation}
For any sequence $C=(c_{n})_{n\geq 0}$,  we denote by $c^{\star k}_{n}=\sum_{i_{1}+\dots+i_{k}=n}c_{i_{1}}\cdots c_{i_{k}}$ 
the $n$th term of the sequence
$C^{\star k}=\overbrace{C\star\cdots\star C}^{\times k}=\left(c^{\star k}_{n}\right)_{n\geq 1}$.
\subsection{Central binomial coefficients}
We consider the generating function,
\begin{equation}
	B(y)=\frac12\left(\frac1{\sqrt{1-4y}}-1\right),
\end{equation}
 appearing in the argument of $\exp$ in the right hand side of equality (\ref{A(y,z)}). It is well known that ${B(y)\over y}$ is the
 ordinary generating function of the  binomial coefficient $\binom{2n+1}{n+1}$. Hence,
 \begin{equation}
 \begin{array}{rcl}
	 \exp\left\{\frac12z\left(\frac1{\sqrt{1-4y}}-1\right)\right\}
	 &=&\exp\left\{zy{B(y)\over y}\right\}
	 =\displaystyle\sum_{k\geq 0}\frac1{k!}\left({B(y)\over y}\right)^{k}z^{k}y^{k}\\
	 &=&\displaystyle\sum_{k,n\geq 0}\frac1{k!}{\binom{2n+1}{n+1}}^{\star k}z^{k}y^{n+k}.
	 \end{array}
 \end{equation}
 Branching this equality in (\ref{A(y,z)}), one obtains
 \begin{equation}
	 \begin{array}{rcl}
		 A(y,z)&=&\displaystyle\sum_{k\geq 0}\frac1{k!}z^{k}y^{k+1}
		 \left(\sum_{n\geq 0} 4^{n}y^{n}\right)\left(\sum_{n\geq 0}{\binom{2n+1}{n+1}}^{\star k}y^{n}\right)\\
		 &=&\displaystyle\sum_{0\leq k<n}\frac1{k!}\left(\sum_{i+j=n-k-1}{4^{i}\binom{2j+1}{j+1}}^{\star k}\right)z^{k}y^{n}.
	 \end{array}
 \end{equation}
 Hence, for a fixed $k$ the sequence $(\tilde\alpha_{k,n})_{n}$ can be written as a convolution.
 \begin{proposition}\label{propalpha}\begin{equation}
	 \tilde\alpha_{k,n}=\frac1{k!}\sum_{i+j=n-k-1}{4^{i}\binom{2j+1}{j+1}}^{\star k}.
	 \end{equation}
 \end{proposition}
 The first values of $\tilde\alpha_{k,n}$ are collected in Table \ref{tablealpha}.
 \begin{example}
	 For instance we have
	 \[{}\begin{array}{rcl}
	 \tilde\alpha_{3,6}&=&\frac1{6!}\left(\binom53^{\star 3}+4\binom32^{\star 3}+16\binom11^{\star 3}\right){}
	 =\frac1{6!}\left(57+4\cdot 9+16\right).\\&=&{109\over 6}.\end{array}
	 \]
 \end{example}
 From Proposition \ref{propalpha}, we observe that
 $ \tilde\alpha_{0,n}=4^{n-1}$ and $\tilde\alpha_{n-1,n}=\frac1{n!}$.
 \begin{table}[H]
 \resizebox{\textwidth}{!}{  
   \begin{tabular}{ |l|| l| l| l| l| l| l| l| l| l| l| l}
   \hline
  {\small $k\setminus n$}&1&  2& 3&  4&    5&   6&     7&           8&        9&   			10\\ \hline
			 0&1& 4& 16& 64& 256& 1024& 4096& 16384&    65536&  262144\\
			 1&0& 1& 7& 38& 187&  874&  3958& 17548& 76627& 330818\\
			 2&0& 0& 1&  10&    69& 406&  2186&  11124&  54445& 259006\\
			 3&0& 0& 0&  1&     13&   109&   748&   4570&   25879&  138917\\
			 4&0& 0& 0&  0&     1&      16&  158&      1240&     8485&      52984 \\
			 5&0& 0& 0&  0&     0&    1&      19& 	216&      1909&	    14471 \\
			 6&0& 0& 0&  0&     0&   0&      1&   	  22&      283& 	   2782 \\
			 7&0& 0& 0&  0&     0&   0&     0&   	  1&      25&    359 \\
			 8&0& 0& 0&  0&     0&   0&     0&  	 0&      1&    28 \\
			 9&0& 0& 0&  0&     0&   0&     0&  	 0&		0&    1  \\  \hline
	 \end{tabular}
	 }
	 \caption{First values of $k!\tilde\alpha_{k,n}.$\label{tablealpha}}
 \end{table}
 The following sequences are registred in OEIS encyclopedia : $\tilde{\alpha}_{1,n}$ is $A000531$ and count the number from area of cyclic polygon of $2n+1$ sides.
  $2!\tilde{\alpha}_{2,n}\quad  (n\geq 3)$ is $A038806$, $i!\tilde{\alpha}_{i,i+3}\quad  (n\geq 3)$ is $A081270$ and $i!\tilde{\alpha}_{i,i+2}\quad  (n\geq 2)$ is $A016777$.  \\
\subsection{Application to Fibonacci numbers and related sequences}
We consider a sequence $A=(a_{n})_{n\geq 0}$ whose generating function is of the form
\begin{equation}
	S_{A}(x)=\sum_{k\geq 0}a_{k}x^{k}={-\gamma\over Q},
\end{equation} 
where $\gamma$ is a complex number and $Q$ is a unitary quadratic polynomial with two distinct roots.
So we have
\begin{equation}
	\sum_{1\leq n}\sum_{0\leq k}a^{\star n}_{k}x^{k}y^{n}=\frac {yS_{A}(x)}{1-yS_{A}(x)}={y\gamma F\over 1-y\gamma F},
\end{equation}
where $F=-\frac1Q$. Hence from equality (\ref{D2F}), one has
\begin{equation}
	\sum_{1\leq n}\sum_{0\leq k}a^{\star n}_{k}x^{k}y^{n}=\mathbb D^{\gamma y}_{x} F.
\end{equation}
Equaling the coefficient of $y^{n}$ in the left hand side and the right hand side, we find
\begin{equation}\label{devS}
	\left(S_{A}(x)\right)^{n}
	=\left(\frac\gamma\Delta\right)^{n-1}\sum_{i=0}^{n-1}\tilde\alpha_{i,n}
	{\left(\diff Qx\right)}^{i}\idiff{S_{A}(x)}x{i}
\end{equation}
where $\Delta$ is the discriminant of $Q$.\\
Now extracting
the coefficient of $x^{n}$ in (\ref{devS}), one finds
\begin{equation}\label{coefxS}
	\begin{array}{rcl}a^{\star n}_{k}&=&\displaystyle
	\left(\frac\gamma\Delta\right)^{n-1}\sum_{i=0}^{n-1}\tilde\alpha_{i,n}\sum_{k_{1}+k_{2}=k}
	\binom i{k_{1}}2^{k_{1}}(-p)^{i+k_{2}-k}(k_{2}+i)_{i}a_{k_{2}+i}
	\\&=&\displaystyle\left(\frac\gamma\Delta\right)^{n-1}\sum_{s=k}^{k+n-1}
	\left(\sum_{i=s-k}^{n-1}\tilde\alpha_{i,n}\binom i{k-s+i}2^{k-s+i}(-p)^{s-k}(s)_{i}\right)a_{s}
\end{array}
\end{equation}
where $Q=x^{2}-px+q$ and $(s)_{i}=s(s-1)\cdots (s-i+1)$ denotes the Pochhammer symbol.\\ \\
Let $\alpha$ and $\beta$ be two distincts non zero complex numbers.
 We consider, as in \cite{ZW}, the generalized Fibonacci numbers defined as $U_{n}={\alpha^{n}-\beta^{n}\over\alpha-\beta}$ and the generalized
 Lucas numbers defined as $V_{n}=\alpha^{n}+\beta^{n}$. Setting $q=\alpha\beta$ and $Q_{k}(x)=q^{-k}-{V_{k}\over q^{k}}x+x^{2}$. One has
 \begin{equation}
	 F_{k}(x):=\sum_{n\geq 0}U_{k(n+1)}x^{n}={{U_{k}\over q^{k}}\over Q_{k}}.
 \end{equation}
Since $\alpha\neq\beta$ the discriminant $\Delta_{k}={V_{k}^{2}-4q^{k}\over q^{2k}}$ of $Q_{k}$ is not zero. So we can apply equality (\ref{devS}) to $F_{k}$
and obtain
\begin{equation}\label{devF}
\begin{array}{rcl}
	\left(F_{k}(x)\right)^{n}
	&=&\displaystyle\left({-U_{k}\over q^{k}}\over\Delta_{k}\right)^{n-1}\sum_{i=0}^{n-1}\tilde\alpha_{i,n}
	{\left(\diff Qx\right)}^{i}\idiff{F_{k}(x)}x{i}\\
	&=&\displaystyle\left(U_{k}\over V_{k}^{2}-4q^{k}\right)^{n-1}\sum_{i=0}^{n-1}(-1)^{n-i-1}q^{(n-i-1)k}\tilde\alpha_{i,n}(V_{k}-2xq^{k})^{i}
	\idiff{F_{k}}x{i}.\end{array}
\end{equation}
\begin{example}\rm
	For instance, for $n=4$ one obtains
	\[{}\begin{array}{l}
	\left(F_{k}\right)^{4}=\left(U_{k}\over V_{k}^{2}-4q^{k}\right)^{3}\left(\frac16{}
	(V_{k}-2xq^{k})^{3}
	\idiff{F_{k}}x{3}-
	5q^{k}{}
	(V_{k}-2xq^{k})^{2}
	\idiff{F_{k}}x{2}\right.\\\left.+38q^{2k}
	(V_{k}-2xq^{k})
	\diff{F_{k}}x-64q^{3k}F\right),\end{array}
	\]
	as expected in \cite{ZW}.\\
	Consider also the bigger example for $n=10$ which gives
	\[{}\begin{array}{l}
	\left(F_{k}\right)^{10}=\left(U_{k}\over V_{k}^{2}-4q^{k}\right)^{9}\left(\frac1{9!}
	(V_{k}-2xq^{k})^{9}
	\idiff{F_{k}}x{9}-
	\frac1{1440}q^{k}{}
	(V_{k}-2xq^{k})^{8}
	\idiff{F_{k}}x{8}\right.\\+\frac{359}{5040}q^{2k}
	(V_{k}-2xq^{k})^{7} \idiff{F_{k}}x7
	-{1391\over 360}q^{3k}
	(V_{k}-2xq^{k})^{6} \idiff{F_{k}}x6
	\\+{14471\over 120}q^{4k}
	(V_{k}-2xq^{k})^{5}\idiff{F_{k}}x5-{6623\over 3}q^{5k}
	(V_{k}-2xq^{k})^{4} \idiff{F_{k}}x4\\
	{138917\over 6}q^{6k}
	(V_{k}-2xq^{k})^{3}\idiff{F_{k}}x3-129503q^{7k}
	(V_{k}-2xq^{k})^{2} \idiff{F_{k}}x2\\
	\left.+330818q^{8k}
	(V_{k}-2xq^{k})\diff{F_{k}}x
	-262144q^{9k}F\right).\\
	\end{array}
	\]
\end{example}

Set $f_{j}^{(k)}=U_{k(j+1)}$. From (\ref{coefxS}), we obtain
\begin{equation}
	\left(f^{(k)}_{j}\right)^{\star n}=\left(U_{k}\over V^{2}_{k}-4q^{k}\right)^{n-1}\sum_{s=0}^{n-1}(-q^{k})^{n-s-1}\beta_{n,s}(j+n)V_{k}^{s}f_{s+j}^{(k)},
\end{equation}
where $\beta_{n,s}(x)=\sum_{i=s}^{n-1}2^{i-s}\binom i{i-s}\tilde\alpha_{i,n}(s-n+x)_{i}$. See in Appendix \ref{AppBeta} for the first values of these polynomials.
Since we have
\begin{equation}
	\left(f_{j}^{(k)}\right)^{\star n}=\sum_{j_{1}+\cdots+j_{n}=j+n}U_{kj_{1}}\cdots U_{kj_{n}},
\end{equation}
we deduce the following result.
\begin{corollary}\label{corU}
	\[{}
	\sum_{j_{1}+\cdots+j_{n}=j}U_{kj_{1}}\cdots U_{kj_{n}}=
	\left(U_{k}\over V^{2}_{k}-4q^{k}\right)^{n-1}\sum_{s=0}^{n-1}(-q^{k})^{n-s-1}\beta_{n,s}(j)V_{k}^{s}U_{k(s+j-n+1)}.
\]
\end{corollary}
\begin{example}\rm
Let us treat the case where $n=6$. Applying Corollary \ref{corU}, one finds
\[{}\begin{array}{l}\displaystyle
\sum_{j_{1}+\cdots+j_{6}=j}U_{kj_{1}}U_{kj_{2}}U_{kj_{3}}U_{kj_{4}}U_{kj_{5}} U_{kj_{6}}=
	\left(U_{k}\over V^{2}_{k}-4q^{k}\right)^{n-1}\left(
	\beta_{6,5}(j)V_{k}^{5}U_{kj}\right.\\\left.-q^{k}\beta_{6,4}(j)V_{k}^{4}U_{k(j-1)}
	+q^{2k}\beta_{6,3}V_{k}^{3}U_{j-2}-q^{3k}\beta_{6,2}V_{k}^{2}U_{j-3}
	+q^{4k}\beta_{6,1}V_{k}U_{j-4}\right.\\\left.-q^{5k}\beta_{6,0}U_{j-5}\right),\end{array}
\]
as in \cite{Mansour} (the values of $\beta_{6,i}(x)$ are in Appendix \ref{AppBeta}).\\
Also consider the bigger example where $n=10$:

\[{}
\begin{array}{l}\displaystyle{}
\sum_{j_{1}+\cdots+j_{10}=j}U_{kj_{1}}\cdots U_{kj_{10}}=
\left(U_{k}\over V^{2}_{k}-4q^{k}\right)^{n-1} \left(
\frac {1}{362880}\, (j-1)_{9} V_{k}^{9}U_{kj}\right.\\
 -\frac{1}{20160}q^{k}(j-2)_{8}(j+4)V_{k}^{8}U_{k(j-1)}+
 \frac1{5040}q^{2k}(j-3)_{7}(2j^2+14j+19)V_{k}^{7}U_{k(j-2)}\\
 -\frac1{1080}q^{3k}(j-4)_{6}(j+3)(2j^{2}+12j+1)V_{k}^{6}U_{k(j-3)}\\
 +\frac1{360}q^{4k}(j-5)_{5}(2j^{4}+20j^{3}+40j^{2}-5j-87)V_{k}^{5}U_{k(j-4)}\\
 -\frac1{180}q^{5k}(j-6)_{4}(2j^{4}+16j^{3}-12j^{2}-176j+75)V_{k}^{4}U_{k(j-5)}\\
 +\frac1{270}q^{6k}(j-7)_{3}\left( 4\,{j}^{6}+36\,{j}^{5}-50\,{j}^{4}-840\,{j}^{3}-404
\,{j}^{2}+3054\,j+1125 \right) V_{k}^{3}U_{k(j-6)}\\

-\frac1{315}q^{7k}(j-8)_{2} \left( j+
1 \right)  \left( 4\,{j}^{6}+24\,{j}^{5}-206\,{j}^{4}-984\,{j}^{3}\right.\\\left.+
2860\,{j}^{2}+7752\,j-11025 \right)V_{k}^{2}U_{k(j-7)} 
\\
+\frac2{315}q^{8k}(j-9)
 \left( {j}^{8}+4\,{j}^{7}-84\,
{j}^{6}-266\,{j}^{5}+1974\,{j}^{4}+4396\,{j}^{3}\right.\\
\left.-12916\,{j}^{2}-15159
\,j+11025 \right) V_{k}U_{k(j-8)}\\-\left.\frac4{2835}q^{9k}
(j^{2}-2^{2})(j^{2}-4^{2})(j^{2}-6^{2})(j^{2}-8^{2})U_{k(j-9)}\right)
\end{array}
\]
\end{example}

\section{Tribonacci, quadrabonacci and beyond}
The technical described in the paper can be generalized for higher level recurrences.
We illustrate our purpose with a level $3$ recurrence. We assume that the generating series
is given by $F=-1/Q$ where $Q=x^{3}+3ax^{2}+3bx+c$ is a unitary cubic polynomial with three distinct 
roots. We consider the binary cubic $P(x,y)=x^{3}+3ax^{2}y+3bxy^{2}+cy^{3}$. Polynomial invariants of binary forms are studied since 
the middle of the Nineteenth century. These invariants are obtained by applying the Cayley Omega process (see e.g. \cite{Olver}).
A transvection is a bilinear operation defined by
\begin{equation}
	(P,Q)^{k}=\left.\left|\begin{array}{cc}{\partial\over\partial x_{1}}&{\partial\over\partial y_{1}}\\ \  \\
	{\partial\over\partial x_{2}}&{\partial\over\partial y_{2}} \end{array}\right|^{k}P(x_{1},y_{1})Q(x_{2},y_{2})\right|_{x_{1}=x_{2}=x\atop{}
	y_{1}=y_{2}=y}.
\end{equation}
For instance, the algebra of {polynomial invariants} of a binary cubic is generated by 
\begin{equation}\label{invariant j}
	j=\frac1{1296}((P,P)^{2},(P,P)^{2})^{2}=6abc-4b^{3}-4a^{3}c+3a^{2}b^{2}-c^{2}.
\end{equation}
From $P(x,y)=y^{3}Q(xy^{-1})$, equality (\ref{invariant j}) gives
\begin{align}
	36Q\pdiff Qx&\left({\partial Q^{2}\over\partial x^{2}}\right)
	+\left(\left({\partial Q^{2}\over\partial x^{2}}\right)^{2}-16\left(\partial Q\over\partial x\right)^{2}\right)\left({\partial Q\over\partial x}\right)^{2}
	\\ \nonumber  &-2\left({\partial Q^{2}\over\partial x^{2}}\right)^{3}Q-108 Q^{2}=108j.
\end{align}
Since ${\partial F\over\partial x}={\partial Q\over\partial x}{F^{2}}$, we obtain
\begin{equation}\begin{array}{rcl}
	F^{2}&=&\frac1{108j}
	\left({}
	\left(\left({\partial^{2} Q\over\partial x^{2}}\right)^{2}-16\left(\partial Q\over\partial x\right)\right)
	\left({\partial Q\over\partial x}\right){}
	{\partial F\over\partial x}\right.\\&&\left.
	+2\left(\left(\partial^{2}Q\over\partial x^{2}\right)^{3}-18
	\left(\partial Q\over\partial x\right){}
	\left(\partial^{2}Q\over\partial x^{2}\right)\right)F+108\right).{}
	\end{array}
\end{equation}
Now, consider the (shifted) {Tribonaci} numbers defined by  {$T_{0}=T_{1}=1$, $T_{2}=2$} and $T_{n}=T_{n-1}+T_{n-2}+T_{n-3}$ for $n>2$.
For this polynomial we have $108j=-176$, ${\partial Q\over\partial x}=3x^{2}+2x+1$, and ${\partial^{2}Q\over\partial x^{2}}=6x+2$. So we obtain,
\begin{equation}\label{TF2}
	F^{2}=\frac1{176}(4(9x^{4}+12x^{3}+16x^{2}+8x+3){\partial F\over\partial x}+8(27(x^{3}+x^{2}+x)+7)F+108).
\end{equation}
In other words, {by extracting the coefficient of $x^{n}$ in the left and right hand sides of (\ref{TF2}), we observe}
\begin{align}
	T_{n}^{\star {2}}=\frac{1}{44}\Big(&3(n+1)T_{n+1}+2(4n+7)T_{n}+2(8n+19)T_{n-1}\\ \nonumber &+2(6n+15)T_{n-2} +(9n+27)T_{n-3}\Big),
\end{align}
for any $n>0$.

In principle, the same strategy can be applied for computing linearization formula for self-convoluted higher Fibonacci-like numbers. Furthermore,
the algebra of {polynomial invariants} is richer from binary quartic forms since it has more than one generator. 
So for any invariant in a Hilbert basis, one can find a family of formulas. The only limit to this process is the computational 
complexity. The formulas obtained not only have an increasing size but  also the algebra of invariants is difficult to describe. 
For binary forms the last works obtained a description for the binary nonic \cite{BP1} and the binary decimic \cite{BP2}. An alternative way should consist
to compute rational invariants \cite{Olver} instead of polynomial invariants. It is well known that the rational invariants are simpler and can  be computed by using the so-called
associated forms \cite{Meyer} which is obtained from the ground form by applying certain substitutions. But now the difficulties is that we have to recover the differential operators.
Anyway, the closed connexion between the classical invariant theory and the formulas involving self-convolute of generalized Fibonacci numbers
need to be investigated and should bring a series of very interesting and deep results.

\appendix
\section{An alternative proof for Theorem \ref{theoalpha}\label{AppProof}}
Assuming (\ref{defalpha}) and applying (\ref{Q2Delta}), one obtains
\begin{equation}\begin{array}{rcl}
	\left({\partial\over\partial x }Q\right)\left({\partial\over\partial x }F^{n}\right)&=&
	nF^{n-1}\left({\partial \over\partial x }F \right) \left({\partial\over\partial x }Q\right)\\
	&=& nF^{n-1}(\Delta F^{2}-4F).
	\end{array}
\end{equation}
So, we have
\begin{equation}
	F^{n+1}=\frac1{n\Delta}\left(\left({\partial\over\partial x}Q\right)\left({\partial \over\partial x }F^{n}\right) +4nF^{n}\right).
\end{equation}
Substituting $F^{n}$ by its expansion (formula (\ref{defalpha})) and identifying 
the coefficient of ${\partial^{i}\over\partial x^{i}}F$ in the left and right hand sides of the resulting 
formula, we find that the family of polynomials $(\alpha_{i,n}(x))_{i,n}$ defined by the recurrence
\begin{equation}\label{recalpha}
\begin{split}
	\alpha_{i,n}=&0\mbox{ for } i<0,\ n<1\mbox{ or }i\geq n,\\
	\alpha_{0,1}=&1,\\
	\alpha_{i,n+1}=&\frac1{{\Delta} n}\left(\left({\partial \over\partial x}Q\right){}
	\left({\partial \over\partial x}\alpha_{i,n}\right)+\left({\partial \over\partial x}Q\right)\alpha_{i-1,n}+4n\alpha_{i,n}\right){}
	\mbox{ in the other cases}
	\end{split}
\end{equation}
satisfies (\ref{defalpha}).\\
Notice that (\ref{defalpha}) define a unique sequence of polynomials. 
A straightforward induction yields 
\begin{equation}\label{alpha0}%
	\alpha_{0,n}=\left(4\over\Delta\right)^{n},
\end{equation}
and
\begin{equation}\label{alphan-1n}
\alpha_{n,n+1}=\frac1{n!}\left(\frac1\Delta{\partial \over\partial x}Q\right)^{n}.
\end{equation}
More generally, one has
\begin{lemma}\label{alphaQ^i}
	The polynomial $\alpha_{i,n}$ is a degree $i$ polynomial which is equal, up to a multiplicative constant, to $\left({\partial \over\partial x}Q\right)^{i}$.
	\end{lemma}
	{\bf Proof}
		We proceed by induction on the pair $(i,n-i)$ lexicographically ordered. 
		If $i=0$ then equality (\ref{alpha0}) implies the result. Suppose now that $i>0$. 
		If $n-i=1$ then the result is a consequence of equality (\ref{alphan-1n}). 
		Now suppose $n-i>1$. From
		$$\alpha_{i,n}=\frac1{(n-1)\Delta}\left(\left({\partial \over\partial x}Q\right){}
	\left({\partial \over\partial x}\alpha_{i,n-1}\right)+\left({\partial \over\partial x}Q\right)\alpha_{i-1,n-1}+4n\alpha_{i,n-1}\right),$${}
	and ${\partial^{2}\over\partial x^{2}}Q=2$ we show that induction hypothesis implies the lemma.
	$\Box$\\ \\
	For simplicity, we set
	\begin{equation}\label{defalphatilde}
		\tilde\alpha_{i,n}={\Delta^{n-1}\over \left({\partial \over\partial x}Q\right)^{i}}\alpha_{i,n}.
	\end{equation}
	According the Lemma \ref{alphaQ^i}, $\tilde\alpha_{i,n}$ is a rational number satisfying the recurrence
	\begin{equation}\label{rectildealpha}
\begin{split}
	\tilde\alpha_{i,n}=&0, \quad \mbox{ for } i<0,\ n<1\mbox{ or }i\geq n,\\
	\tilde\alpha_{0,1}=&1,\\
	\tilde\alpha_{i,n+1}=&\frac{2i}{n}\tilde\alpha_{i,n}+\frac1n\tilde\alpha_{i-1,n}+4\tilde\alpha_{i,n},\quad {}
	\mbox{ in the other cases }.
\end{split}
	\end{equation}
	Consider the double generating function  
	{ \begin{equation}
	A(y,z)=\sum_{0 \leq i<n<\infty}\tilde\alpha_{i,n}z^{i}y^{n},\end{equation}}
	Equivalently to (\ref{rectildealpha}), the series $A(y,z)$ is the unique solution of
	{
	\begin{equation}
		A(y,z)=y+2yz{\partial\over\partial z}\int_{0}^{y}{A(t,z)\over t}dt+zy\int_{0}^{y}{A(t,z)\over t}dt+4yA(y,z).
	\end{equation}}
	It is  easier to manipulate this equation under the form
	\begin{equation}\label{equadiff}
		(1-4y){\partial^2\over \partial y^{2}}A(y,z)=z\left(1+2{\partial\over\partial z}\right){A(y,z)\over y}
		+(z+8){\partial\over\partial y}
		A(y,z)+2z{\partial^{2}\over\partial y\partial z}A(y,z),
	\end{equation}
	with the initial conditions $A(0,z)=0$, $\left({\partial \over\partial y}A\right)(0,z)=1$, and $A(y,0)={y\over 1-4y}$.
	We check that 
	\begin{equation}\label{soluce}
		A(y,z)={y\over1-4y}\exp\left\{\frac12z\left(\frac1{\sqrt{1-4y}}-1\right)\right\}.
	\end{equation}
	is the unique solution of (\ref{equadiff}).\\
	Hence,
	\begin{equation}
		\begin{array}{rcl}
			\sum_{n\geq 1}F^{n}y^{n}&=&\Delta\left(\sum_{n\geq1}\sum_{i=0}^{n-1}
			\tilde\alpha_{i,n}(x)\left({\partial Q\over\partial x}\right)^{i}{\partial^{i}\over\partial x^{i}}\left(y\over \Delta\right)^{n}\right)F\\
			&=& \Delta A\left({ y\over\Delta},X{\partial\over\partial x}\right)|_{X={\partial Q\over\partial x}}.F\\*
			&=& \frac y{ 1-4{ y\over\Delta}}
			F\left(x+\frac{1}{2}
			{\partial Q\over\partial x}
			\left(\frac1{\sqrt{1-4{y\over\Delta}}}-1\right)\right).
		\end{array}
	\end{equation}

\section{First values of $\beta_{n,p}(x)$ \label{AppBeta}}
In this section we compute the first values for the polynomial $\beta_{n,p}(x)$ involved in Corollary \ref{corU}.
n=1:\\ 
\begin{equation*}
\begin{split}
\beta_{1,0}(x)&=1.\\
\end{split} 
\end{equation*}
n=2:\\ 
\begin{equation*} 
\begin{split}
\beta_{2,0}(x)&=2x,\\ 
		   \beta_{2,1}(x)&=x-1. \\
\end{split} 
\end{equation*}
n=3:\\ 
\begin{equation*} 
\begin{split}
\beta_{3,0}(x)&=2(x+1)(x-1),\\
	       \beta_{3,1}(x)&=(2x+1)(x-2),\\
	       \beta_{3,2}(x)&=\frac12(x-1)(x-2).\\
\end{split}
\end{equation*}
n=4:\\ 
\begin{equation*} 
\begin{split}
\beta_{4,0}(x)&=\frac43x(x+2)(x-2),\\
    	   \beta_{4,1}(x)&=2\, \left( {x}^{2}+x-1 \right) \left( x -3\right),\\
	{}\beta_{4,2}(x)&= \left( x+1 \right)\left( x-2 \right)  \left( x-3 \right)  ,\\{}
	       \beta_{4,3}(x)&=\frac16\, \left( x-1 \right)  \left( x-2 \right)  \left( x-3 \right).\\
\end{split} 
\end{equation*}
n=5:\\ 
\begin{equation*}
\begin{split}
	       \beta_{5,0}(x)&=\frac23(x+3)(x+1)(x-1)(x-3),\\
	       \beta_{5,1}(x)&=\frac13\, \left( 2\,x+1 \right)\left( 2\,{x}^{2}+2\,
x-9 \right)  \left( x-4 \right)  ,\\ 
			\beta_{5,2}(x)&=\frac12\,  \left( 2\,{x}^{2}+4\,x-
1 \right)\left( x-3 \right) \left( x-4 \right),  \\
			\beta_{5,3}(x)&=\frac16\, \left( 2\,x+3 \right)  \left( x-2 \right)  \left( x-3 \right) 
 \left( x-4 \right),\\
			\beta_{5,4}(x)&=\frac1{24}\, \left( x-1 \right)  \left( x-2 \right)  \left( x-3 \right) 
 \left( x-4 \right). \\  
\end{split} 
\end{equation*}
 n=6:\\ 
\begin{equation*}
\begin{split}
			\beta_{6,0}(x)&=\frac4{15}x(x+2)(x+4)(x-2)(x-4),\\
			\beta_{6,1}(x)&=\frac 23\,  \left( {x}^{4}+2\,{x}^{3}-10\,{x}^{2}-11\,x+
9 \right)\left( x-5 \right)  ,\\
			\beta_{6,2}(x)&=\frac13\left( x+1 \right) \left( 2\,{x}^{2}+4\,x-9 \right)  \left( x-4 \right)  \left( x-5 \right)  
,\\
\end{split} 
 \end{equation*}
\begin{equation*} 
\begin{split}
			\beta_{6,3}(x)&=\frac16\,  \left( 2\,{x}^{2}+6\,x+1 \right)\left( x-3 \right)  \left( x-4 \right)  \left( x-5 \right) 
,\\
			\beta_{6,4}(x)&=\frac1{12}\,\left( x+2 \right) \left( x-2 \right)  \left( x-3 \right)  \left( x-4 \right)  \left( x-5 \right)  , \\
			\beta_{6,5}(x)&={\frac {1}{120}}\, \left( x-1 \right)  \left( x-2 \right)  \left( x-3
 \right)  \left( x-4 \right)  \left( x-5 \right) .\\ \end{split} 
 \end{equation*}
 n=7:\\
 \begin{equation*}
\begin{split}
 			\beta_{7,0}(x)&=\frac4{45}(x+1)(x+3)(x+5)(x-1)(x-3)(x-5),\\
			\beta_{7,1}(x)&= \frac2{15}\, \left( 2\,x+1 \right) \left( {x}^{4}+2\,{
x}^{3}-21\,{x}^{2}-22\,x+75 \right) \left( x-6 \right)  ,\\
			\beta_{7,2}(x)&=\frac13\,\left( {x}^{4}+4\,{x}^{
3}-7\,{x}^{2}-22\,x+3 \right) \left( x-5 \right)  \left( x-6 \right)   ,\\
			\beta_{7,3}(x)&=\frac19\,\left( 2\,x+3 \right)  \left( x+4 \right)  \left( x-1 \right)  \left( x-4 \right)  \left( x-5 \right) 
 \left( x-6 \right)   ,\\
			\beta_{7,4}(x)&=\frac1{12}\, \left( {x}^{2}+4\,x+2 \right) \left( x-3 \right)  \left( x-4 \right)  \left( x-5 \right) 
 \left( x-6 \right)  ,\\
\beta_{7,5}(x)&={\frac {1}{120}}\, \left( 2\,x+5 \right) \left( x-2 \right)  \left( x-3 \right)  \left( x-4 \right)  \left( x-5 \right)  \left( x-6 \right) ,\\
			\beta_{7,6}(x)&= {\frac {1}{720}}\, \left( x-1 \right)  \left( x-2 \right)  \left( x-3
 \right)  \left( x-4 \right)  \left( x-5 \right)  \left( x-6 \right).\\
\end{split} 
\end{equation*}

\vspace{0.5 cm}
It is easy to prove that $\beta_{n,k}$ is a degree  $n-1$ polynomial admitting 
$(x-n+k)_{k}$ as a factor but the polynomials $\beta_{n,k}(x-n+k)_{k}^{-1}$ remains to be
properly investigated. For instance, numerical evidence suggests an interesting closed form for  
$\beta_{n,0}$.\\ \\

\noindent{\bf Acknowledgment:}
The paper was partially supported by the project PHC MAGHREB ITHEM 14MDU929M and GRR project MOUSTIC.

\end{document}